\documentclass[12pt]{amsart}
\usepackage{amssymb,latexsym, amsmath, amsxtra}
\usepackage{graphicx}
\begin{document}

\title{Moments of zeta and correlations of divisor-sums: III}
 
\author{Brian Conrey}
\address{American Institute of Mathematics, 360 Portage Ave, Palo Alto, CA 94306, USA and School of Mathematics, University of Bristol, Bristol BS8 1TW, UK}
\email{conrey@aimath.org}
\author{Jonathan P. Keating}
\address{School of Mathematics, University of Bristol, Bristol BS8 1TW, UK}
\email{j.p.keating@bristol.ac.uk}

\thanks{We gratefully acknowledge support under EPSRC Programme Grant EP/K034383/1
LMF: L-Functions and Modular Forms.  Research of the first author was also supported by the American Institute of Mathematics and by a grant from the National Science Foundation. JPK is grateful for the following additional support: a grant from the Leverhulme Trust, a Royal Society Wolfson Research 
Merit Award, a Royal Society Leverhulme Senior Research Fellowship, and a grant from the Air Force Office of Scientific Research, Air Force Material Command, 
USAF (number FA8655-10-1-3088). He is also pleased to thank the American Institute of Mathematics for hospitality during a visit where this work started.}

\date{\today}

\begin{abstract} In this series 
we examine the calculation of the $2k$th  moment and shifted moments of the Riemann zeta-function on the critical line using long Dirichlet 
polynomials and divisor correlations.   The present paper is concerned with the precise input of the conjectural formula for the 
classical shifted convolution problem for
  divisor sums  so as to obtain all of the lower order terms in the asymptotic formula for the mean square along
 $[T,2T]$ of a Dirichlet polynomial of length 
up to $T^2$ with divisor functions as coefficients.
\end{abstract}

\maketitle

\section{Introduction}
This paper is part 3 of a sequence of papers devoted to understanding how to conjecture all of 
the integral moments of the Riemann zeta-function from a number theoretic perspective. The method is to approximate
$\zeta(s)^k$ by a long Dirichlet polynomial and then compute the mean square of the Dirichlet polynomial (c.f.~[GG]). 
There will be many off-diagonal terms and it is the care of these that is the concern of these papers. In
particular it is necessary to treat the off-diagonal terms by a method invented by Bogomolny and Keating [BK1, BK2]. Our perspective on
this method is that it is most properly viewed as a multi-dimensional Hardy-Littlewood circle method. 

In parts 1 and 2 [CK1, CK2] we considered the 
second and fourth moments of zeta in this new light. In this paper we embark on the higher moments. Here we only consider the 
classical shifted convolution problem but we try to solve it precisely in a way that exhibits all of the 
main terms of the expected formula. This requires some care! In the next paper we will introduce  the new terms for the higher moments
(much as we did for the fourth moment in paper 2 [CK2]). 

One way to think of this paper is that it is a more precise version
of [CG] in that we obtain all of the main terms in the asymptotic formula. Our treatment is not rigorous; in particular we 
conjecture the shape of the fundamental shifted convolution at a critical juncture. This is to be expected since for example
no one knows how to evaluate
$$\sum_{n\le x } \tau_3(n)\tau_3(n+1)$$
asymptotically in a rigorous way.

The  formula we obtain is in complete agreement with all of the main terms  predicted by the recipe of [CFKRS] (and in particular, with the leading order term conjectured in [KS]). 
 
\section{Shifted moments}
Some of the underlying mechanism of moments becomes a little clearer if we introduce shifts. (The initial work is possibly harder
but the payoff makes it worthwhile.) Basically we are interested in the moment
\begin{eqnarray*} \int_0^T \zeta(s+\alpha_1)\dots \zeta(s+\alpha_k)\zeta(1-s+\beta_1)\dots \zeta(1-s+\beta_k) ~dt
\end{eqnarray*}
where $s=1/2+it$ and we think of the shifts as being small complex numbers of size $\ll 1/\log T$.
Now
$$\zeta(s+\alpha_1)\dots \zeta(s+\alpha_k)=\sum_{m_1,\dots,m_k} \frac{1}{m_1^{s+\alpha_1}\dots m_k^{s+\alpha_k}}=\sum_{m=1}^\infty
\frac{\tau_{\alpha_1,\dots , \alpha_k}(m)}{m^s}$$
where 
$$ \tau_{\alpha_1,\dots , \alpha_k}(m)=\sum_{m_1\cdot m_2 \dots  m_k=m} m_1^{-\alpha_1}\dots m_k^{-\alpha_k}.$$
(We have here used $\tau$ for the divisor function rather than $d$.) Let 
$$D_{\alpha_1,\dots,\alpha_k}(s;X)=\sum_{n\le X} \frac{\tau_{\alpha_1,\dots,\alpha_k}(n)}{n^s}.$$
More succinctly  
if $A=\{\alpha_1,\dots, \alpha_k\}$ we let
$$\tau_A(m)=\tau_{\alpha_1,\dots,\alpha_k}(m)$$
and
$$D_A(s) =\sum_{m=1}^\infty \frac{\tau_A(m)}{m^s}=\prod_{\alpha\in A} \zeta(s+\alpha);$$
for the Dirichlet polynomial approximation to this we use the notation
$$D_A(s;X) =\sum_{m\le X} \frac{\tau_A(m)}{m^s}.$$
For a set $A$ it is convenient to designate this set translated by $w$ with the notation $A_w$; i.e.
 $$A_w=\{w+\alpha_1,\dots, w+\alpha_k\}.$$

The moment we are interested in is 
\begin{eqnarray}
\label{eqn:moment} 
R^\psi_{A,B}(T)=\int_0^\infty \psi\left(\frac t T \right)  \prod_{\alpha\in A}\zeta(s+\alpha)\prod_{\beta\in B} \zeta(1-s+\beta) ~dt
\end{eqnarray}
where $\psi$ is a smooth function with compact support, say $\psi \in C^\infty[1,2]$.
The recipe [CFKRS] tells us how we  expect this moment will behave, namely 
$$ R^\psi_{A,B}(T)=T\int_0^\infty \psi(t)\sum_{U\subset A, V\subset B\atop |U|=|V|}
 \left(\frac{tT}{2\pi}\right)^{-\sum_{\hat{\alpha}\in U\atop
\hat{\beta}\in V}(\hat{\alpha}+\hat{\beta}) }\mathcal A Z(A-U+V^-,B-V+U^-)~dt+o(T)
$$
where 
$$Z(A,B):=\prod_{\alpha\in A , \beta\in B} \zeta(\alpha+\beta)$$
and $\mathcal A (A,B)$ is a product over primes that converges nicely
in the domains under consideration (see below).  Also we have used an unconventional notation here; 
 by $A-U+V^-$ we mean the following: start with the set $A$ and remove the elements of $U$ and then include the negatives of 
the elements of $V$. We think of the process as ``swapping" equal numbers of elements between  $A$ and $B$; when elements are removed from $A$
and put into $B$ they first get multiplied by $-1$. We keep track of these swaps with our equal-sized subsets $U$ and $V$  of $A$ and $B$;
and when we refer to the ``number of swaps'' in a term we mean the cardinality $|U|$ of $U$ (or, since they are of equal size, of $V$). 

The Euler product $\mathcal A$ is given by
 \begin{eqnarray*}\label{eq:Azeta}
  \mathcal{A} (A,B)=\prod_p Z_{p}(A,B)\int_0^1\mathcal{A}_{p,\theta}(A,B)~d\theta,
  \end{eqnarray*}
  where
$z_p(x):=(1-p^{-x})^{-1}$, $Z_p(A,B)=\prod_{\alpha\in
A\atop\beta\in B} z_p(1+\alpha+\beta)^{-1}$  
  and
\begin{eqnarray*}
\label{eq:Aptheta} \mathcal{A}_{p,\theta}(A,B):=  
\prod_{\alpha\in A} z_{p,-\theta}(\frac 12 +\alpha)
\prod_{\beta\in B}z_{p,\theta}(\frac 12 +\beta)
\end{eqnarray*}
with $z_{p,\theta}(x):=(1-e(\theta)p^{-x})^{-1}$.

The technique we are developing here is to approach our  moment problem (\ref{eqn:moment}) through long Dirichlet polynomials, i.e.
 we consider  
\begin{eqnarray} \label{eqn:basicI}  \nonumber
I^\psi_{A,B}(T;X):&=&\int_0^\infty \psi\left(\frac t T\right) D_{A}(s;X)D_{B}(1-s;X) ~dt\\
&=& T\sum_{m,n\le X }\frac{\tau_A(m)\tau_B(n)\hat{\psi}\left(\frac{T}{2\pi}\log \frac m n \right)}{\sqrt{mn}}
\end{eqnarray}
for various ranges of $X$. We can use the recipe of [CFKRS] to conjecture a formula for $I^\psi$. We start with  
$$ D_{A}(s;X)=\frac{1}{2\pi i}
\int_{w} \frac{X^{w}}{w}D_{A_w}(s)~dw.
$$
Thus,
$$I^\psi_{A,B}(T;X)= \frac{1}{(2\pi i)^2}
\iint_{z,w}\frac{X^{z+w}}{zw} R^\psi_{A_{w},B_{z}}(T)~dw~dz.
$$
 We insert the conjecture above from the recipe and expect that
\begin{eqnarray*}
 I^\psi_{A,B}(T;X) &=&  T\int_0^\infty \psi(t)\frac{1}{(2\pi i)^2}
\iint_{z,w}\frac{X^{z+w}}{zw}\sum_{U\subset A, V\subset B\atop |U|=|V|}
 \left(\frac{tT}{2\pi}\right)^{-\sum_{\hat{\alpha}\in U\atop
\hat{\beta}\in V}(\hat{\alpha}+w+\hat{\beta}+z) }\\
&&\qquad \qquad \times \mathcal A Z(A_w-U_w+V^-_z,B_z-V_z+U^-_w)~dw~dz~dt +o(T).
\end{eqnarray*}
We have done a little simplification here: instead of writing $U\subset A_w$ we have written
$U\subset A$ and changed the exponent of $(tT/2\pi)$ accordingly.

Notice that there is a factor $(X/T^{|U|})^{w+z}$ here. As mentioned above we refer to  $|U|$ as the number of ``swaps" in the recipe, and now we see
more clearly the role it plays; in the terms above for which  $X<T^{|U|}$   we move the path of 
integration in $w$ or $z$ to $+\infty$ so that the factor  $(X/T^{|U|})^{w+z}\to 0$ and the contribution of such a term is 0.
Thus, the size of $X$ determines how many ``swaps'' we must keep track of. For example, 
if $X<T$, then we only need to
keep the terms with $|U|=|V|=0$, i.e. no swaps. We then have
\begin{eqnarray} \label{eqn:I}
I^\psi_{A,B}(T;X)=T\hat{\psi}(0) \frac{1}{(2\pi i)^2}\iint_{\Re z=2\atop\Re w=2}
\frac{X^{w+z}}{wz} \mathcal A Z(A_w,B_z)~dw ~dz +o(T).
\end{eqnarray}
We let $s=z+w$ and have
$$I^\psi_{A,B}(T;X)=T\hat{\psi}(0)\frac{1}{(2\pi i)^2}\iint_{\Re s=4\atop\Re w=2}
\frac{X^{s}}{w(s-w)} \mathcal A Z(A_s,B)~dw ~ds +o(T);$$
(here we used an identity $\mathcal A Z(A_w,B_z)=\mathcal A Z(A_{w+z},B)$; this is obvious for the $Z$ 
factor and less obvious for the $\mathcal A$ factor).
We move the $w$ integral to the left towards $\Re w = -\infty$ and evaluate that integral as the residue at $w=0$.
Thus,
$$I^\psi_{A,B}(T;X)=T\hat{\psi}(0)\frac{1}{2\pi i}\int_{\Re s=4}
\frac{X^{s}}{s}\mathcal A  Z(A_s,B) ~ds +o(T).$$
Since 
\begin{eqnarray*}
\sum_{n=1}^\infty \frac{\tau_A(n)\tau_B(n)}{n^s}=\mathcal A  Z(A_s,B)
\end{eqnarray*}
it follows that
$$I^\psi_{A,B}(T;X)=T\hat{\psi}(0)\sum_{n\le X} \frac{\tau_A(n)\tau_B(n)}{n} +o(T).$$ 
This is exactly the same as what we have referred to previously [CK1, CK2] as the diagonal term.

\section{Type I off-diagonals with shifts}

What happens if $T\ll X \ll T^2$?  From the recipe we add on the terms with one swap, i.e. $|U|=|V|= 1.$
These terms are
\begin{eqnarray*} && \label{eqn:I}
 T\int_0^\infty \psi(t) \sum_{\hat\alpha\in A \atop \hat \beta\in  B }
\left(\frac{Tt}{2\pi}\right)^{- \hat{\alpha}-
\hat{\beta} } 
 \frac{1}{(2\pi i)^2}\iint_{\Re z=2\atop\Re w=2}
\frac{X^{w+z}}{wz} \\
&&\qquad \qquad \mathcal A Z\big((A_w-\{\hat \alpha+w\}) \cup \{-\hat\beta-z\},(B_z-\{\hat \beta+z\})\cup 
\{-\hat \alpha-w\}\big) 
~dw ~dz ~dt +o(T).
\end{eqnarray*}
Let $A'=A-\{\hat\alpha\}$ and $B'=B-\{\hat \beta\}$. Then  
\begin{eqnarray*}&&
  Z\big((A_w-\{\hat \alpha+w\}) \cup \{-\hat\beta-z\},(B_z-\{\hat \beta+z\})\cup 
\{-\hat \alpha-w\}\big) \\&&\qquad 
=  Z(A'_w,B'_z)   Z(A'_w,\{-\hat \alpha -w\})   Z(\{-\hat\beta-z\},B'_z)   Z(\{-\hat \beta-z\},\{-\hat \alpha-w\})\\
&&\qquad =   Z(A'_{z+w},B')  Z(A',\{-\hat \alpha\})  Z(\{-\hat \beta\},B') \zeta(1-\hat\alpha-\hat\beta -w-z).
\end{eqnarray*}
Letting $s=w+z$ and integrating the $w$-integral, we obtain 
\begin{eqnarray} \label{eqn:oneswap}&&
T\int_0^\infty \psi(t) \sum_{\hat\alpha\in A\atop \hat \beta \in B}
\left(\frac {Tt}{2\pi}\right)^{-\hat \alpha -\hat \beta}  Z(A',\{-\hat \alpha\})  Z(\{-\hat \beta\},B')\\
&& \qquad \qquad \times \frac{1}{2\pi i}
\int_{\Re s=4}
\frac{\left(\frac {2\pi X}{Tt}\right)^s}{s} \mathcal A(A'\cup\{-\hat\beta-s\},B'_s\cup\{-\hat\alpha\})
 Z(A'_{s},B')  \zeta(1-\hat\alpha-\hat\beta -s)~ds. \nonumber
\end{eqnarray}
 
How does such a term appear on the coefficient correlation side of things?
The correlation sum 
$$D_{A,B}(u,h):=\sum_{n\le u}  \tau_A(n) \tau_B(n+h)$$
when averaged over $h$ turns out to lead to just such an expression.
  We conjecture that (see [CG] and [DFI])
$$
D_{A,B}(u,h)=m_{A,B}(u,h)+O(u^{1/2+\epsilon}) 
$$
uniformly for $1\le h \le u^{1-\epsilon}$, where $m_{A,B}(u,h)$ is
a smooth function of $u$ whose derivative  is  
$$
m_{A,B}'(u,h)=\sum_{d\mid h} \frac{f_{A,B}(u,d)}{d} ,  
$$
where
$$
f_{A,B}(u,d)=\sum_{q=1}^\infty \frac{\mu(q)}{q^2}P_A(u,qd)P_B(u+h,qd) ,  
$$
in which $P_A(u,q)$ is the average of $\sum_{n\le u}\tau_A(n) e(n/q)$, i.e.
$$
P_A(u,q)=\frac{1}{2\pi i}\int_{|s|=1/8}
\prod_{\alpha\in A}\zeta(s+1+\alpha)G_A(s+1,q)\left(\frac{u}{q}\right)^s~ds\ ,  
$$
with
$$
G_A(s,q)=\sum_{d\mid q}\frac{\mu(d)}{\phi(d)}d^s\sum_{e\mid d}
\frac{\mu(e)}{e^s}g_A(s,qe/d)\ , 
$$
and, if $q=\prod_pp^{q_p}$,
$$
g_A(s,q)=\prod_{p\mid
q}\left(\prod_{\alpha\in A}\left(1-p^{-s-\alpha}\right)\sum_{j=0}^\infty
\frac{\tau_A\left(p^{j+q_p}\right)}{p^{js}}\right)\ .  
$$
 Thus,
$$
  P_A(u,q)=\sum_{\hat \alpha \in A} 
G_A(1-\hat \alpha ,q)\left(\frac{u}{q}\right)^{-\hat\alpha} \prod_{\alpha\in A'}\zeta(1-\hat \alpha +\alpha). 
$$ 
So, looking back to (\ref{eqn:basicI}),  we have to consider 
\begin{eqnarray*}
\sum_{m=n+h\atop T\le m\le X}\frac{\tau_A(m)\tau_B(n)}{m}\hat{\psi}\left(\frac{Th}{2\pi m}\right) &\sim&
\sum_h \int_T^X \langle \tau_A(m)\tau_B(n)\rangle _{m\sim u}^* \hat{\psi}\left(\frac{Th}{2\pi u}\right)\frac{du}{u}
\end{eqnarray*}
where $*$ indicates the condition $m=n+h$.
We replace $\langle \tau_A(m)\tau_B(n)\rangle _{m\sim u}^*$ with 
$$\sum_{q=1}^\infty \frac{\mu(q)}{q^2} \sum_{d\mid h} \frac{1}{d} P_A(u,dq)P_B(u,dq)$$
 and then  switch the sums over $h$ and 
$d$. Thus, the above is 
\begin{eqnarray*}&&
\sum_{\hat{\alpha}\in A\atop \hat{\beta}\in B}
\prod_{\alpha\ne \hat{\alpha}}\zeta(1+\alpha-\hat \alpha)\prod_{\beta\ne \hat{\beta}}\zeta(1+\beta-\hat \beta)
\sum_{q,d,h}  \frac{\mu(q)G_A(1-\hat \alpha, qd) G_B(1-\hat \beta ,qd)}{d^{1-\hat \alpha -\hat \beta}
q^{2-\hat\alpha-\hat \beta}}
\\
&& \qquad \qquad \qquad \times 
 \int_T^X u^{-\hat \alpha -\hat \beta}\hat\psi \left(\frac{Thd}{2\pi u}\right)\frac{du}{u}.
\end{eqnarray*}
We make the change of variable $v=\frac{Thd}{2\pi u}$ and bring the sum over $h$ to the inside; $u<X$ implies that 
$$ hd< \frac{2\pi Xv}{T}.$$ 
Thus, we have 
\begin{eqnarray*}&&
\sum_{\hat{\alpha}\in A\atop \hat{\beta}\in B}
\prod_{\alpha\ne \hat{\alpha}}\zeta(1+\alpha-\hat \alpha)\prod_{\beta\ne \hat{\beta}}\zeta(1+\beta-\hat \beta)\left(\frac
{T}{2\pi}\right)^{-\hat \alpha -\hat \beta}
\int_v \hat{\psi}(v) v^{\hat \alpha +\hat \beta}\\
&& \qquad \qquad 
\times   \sum_{hd< \frac{2\pi Xv}{T}} 
\frac{\mu(q) G_A(1-\hat \alpha, qd) G_B(1-\hat \beta ,qd)}{ q^{2-\hat\alpha-\hat \beta}
 h^{\hat \alpha +\hat \beta}d }\frac{dv}{v}.
\end{eqnarray*}
We use Perron's formula to write the sum over $d,q$ and $h$ as
\begin{eqnarray}  \label{eqn:Perron} &&
\frac{1}{2\pi i}\int_{(2)} \sum_{h,d,q} 
\frac{\mu(q) G_A(1-\hat \alpha, qd) G_B(1-\hat \beta ,qd)}{q^{2-\hat\alpha-\hat \beta}
h^{s+\hat \alpha+ \hat \beta}d^{1+s} }\frac{(2\pi Xv/T)^s}{s}~ds.
\end{eqnarray}
Recall that
$$
G_A(s,q)=\sum_{d\mid q}\frac{\mu(d)}{\phi(d)}d^s\sum_{e\mid d}
\frac{\mu(e)}{e^s}g_A(s,qe/d)\ ,
$$
and
$$
g_A(s,q)=\prod_{p\mid
q}\left(\prod_{\alpha\in A}\left(1-p^{-s-\alpha}\right)\sum_{j=0}^\infty
\frac{\tau_A\left(p^{j+q_p}\right)}{p^{js}}\right)  
$$
for $q=\prod_pp^{q_p}$,
so that
\begin{eqnarray*}
G_A(1-\hat \alpha , p)&=&g_A(1-\hat \alpha,p) -\frac{p^{1-\hat \alpha}}{p-1}+\frac{g_A(1-\hat\alpha,p)}{p-1}\\
&=& \frac{p}{p-1}\left( g_A(1-\hat \alpha,p)-p^{-\hat \alpha}\right)\\
&=& \prod_{\alpha\in A}(1-p^{-1-\alpha+\hat\alpha})\sum_{j=0}^\infty \frac{\tau_A(p^{j+1})}{p^{j(1-\hat\alpha)}}-p^{-\hat\alpha}+O(1/p)\\
&=&\tau_A(p)-p^{-\hat\alpha}+O(1/p) = \tau_{A'}(p)+O(1/p).
\end{eqnarray*}
Thus,
\begin{eqnarray*}
\sum_{d,q=1}^\infty \frac{\mu(q)G_A(1-\hat \alpha, qd) G_B(1-\hat \beta ,qd)}{d^{1+s}q^{2-\hat\alpha-\hat\beta}} &=&
\prod_p \sum_{d,q=0}^\infty \frac{\mu(p^q)G_A(1-\hat \alpha, p^{d+q})G_B(1-\hat \beta,p^{d+q})}{p^{d+dz+2q+qw}}\\
 &=&\mathcal A_{A,B,\hat \alpha ,\hat \beta}(s)  \prod_{a\in A'\atop b\in B'}\zeta(1+a+b+s)
\end{eqnarray*}
where 
$$\mathcal A_{A,B,\hat \alpha ,\hat \beta}(s)
= \prod_p \left(\prod_{\alpha\in A'\atop \beta
\in B'}(1-p^{-1-\alpha-\beta-s})\right)\sum_{d,q=0}^\infty \frac{\mu(p^q)G_A(1-\hat \alpha, p^{d+q})G_B(1-\hat \beta,p^{d+q})}
{p^{d+ds+q(2-\hat\alpha-\hat\beta)}}$$ 
is an Euler product that is absolutely convergent when the real parts of $z$ and $w$ are near 0.
Thus, (\ref{eqn:Perron}) becomes
\begin{eqnarray*}
\\
&&\qquad \qquad 
\frac{1}{2\pi i}\int_{(2)} \zeta(s+\hat \alpha+ \hat \beta)\prod_{a\in A'\atop b\in B'}\zeta(1+s+a+b)
\mathcal A_{A,B,\hat\alpha,\hat \beta}(s)
\frac{\left(\frac{2\pi Xv}{T}\right)^s}{s}~ds.
\end{eqnarray*}

Now
$$\int_v \hat{\psi}(v) v^{\hat \alpha +\hat \beta+s}~\frac{dv}{v}=
\int_0^\infty \psi(t) t^{-\hat \alpha -\hat \beta-s}\chi(1-s-\hat \alpha -\hat \beta).  $$
Thus, the above is 
\begin{eqnarray*}
&&
\sum_{\hat{\alpha}\in A\atop \hat{\beta}\in B}
\prod_{\alpha\ne \hat{\alpha}}\zeta(1+\alpha-\hat \alpha)\prod_{\beta\ne \hat{\beta}}\zeta(1+\beta-\hat \beta)
\left(\frac{T}{2\pi} \right)^{-\hat \alpha -\hat \beta}
\int_0^\infty  \psi(t) t^{-\hat \alpha -\hat \beta} \\
&&\qquad \qquad \times 
\frac{1}{2\pi i}\int_{(2)} 
 \zeta(1-s-\hat \alpha- \hat \beta)
Z(A'_s,B') 
\mathcal A_{A,B,\hat\alpha,\hat\beta}(s)\frac{\left(\frac {2\pi X}{tT}\right)^s}{s}~ds ~dt.
\end{eqnarray*}

This should be compared with (\ref{eqn:oneswap}):
\begin{eqnarray*}  &&
T\int_0^\infty \psi(t) \sum_{\hat\alpha\in A\atop \hat \beta \in B}
\left(\frac {Tt}{2\pi}\right)^{-\hat \alpha -\hat \beta}  Z(A',\{-\hat \alpha\})  Z(\{-\hat \beta\},B')\\
&& \qquad \qquad \times \frac{1}{2\pi i}
\int_{\Re s=4}
\frac{\left(\frac {2\pi X}{Tt}\right)^s}{s} \mathcal A(A'\cup\{-\hat\beta-s\},B'_s\cup\{-\hat\alpha\})
 Z(A'_{s},B')  \zeta(1-\hat\alpha-\hat\beta -s)~ds. \nonumber
\end{eqnarray*}
 These are identical provided that 
 $$
 \mathcal A_{A,B,\hat\alpha,\hat\beta}(s )
=
\mathcal A(A'\cup\{-\hat\beta-s\},B'_s\cup\{-\hat\alpha\}); 
$$
so it just remains to prove this identity. 

\section{The Euler products}
 To prove the identity we follow the method of [CG].
  We compare the $p$-factor of each Euler product
$$\mathcal A_{A,B,\hat \alpha ,\hat \beta}^p(s )
:=    \prod_{\alpha\in A'\atop \beta
\in B'}(1-p^{-1- \alpha- \beta-s}) \sum_{d,q=0}^\infty \frac{\mu(p^q)G_A(1-\hat \alpha, p^{d+q})G_B(1-\hat \beta,p^{d+q})}
{p^{d+ds+q(2-\hat\alpha-\hat \beta)}}$$ 
and 
\begin{eqnarray*}&&
\mathcal A^p(A'\cup\{-\hat\beta-s\},B'_s\cup\{-\hat\alpha\}):= A_p(A'_{s},B')  A_p(A',\{-\hat \alpha\}) A_p(\{-\hat \beta\},B') 
(1-p^{-1+\hat\alpha+\hat\beta +s})\\
&& \qquad \qquad \qquad \qquad \qquad \qquad \qquad \qquad \qquad \times\int_0^1 \mathcal A_{p,\theta}(A'\cup\{-\hat\beta-s\},B'_s\cup\{-\hat\alpha\})
~d\theta
\end{eqnarray*}
Now the integral over $\theta$ is
\begin{eqnarray*}
\sum_{j=0}^\infty \frac{\tau_{A'\cup\{-\hat\beta-s\}}(p^j)\tau_{B'_s\cup\{-\hat\alpha\}}(p^j)}{p^j}.
\end{eqnarray*}
Thus, cancelling the factor $\mathcal A_p(A',B'_s)$ from both sides and replacing all of the $\beta+s\in B_s$ by $\beta$ (i.e. taking $s=0$), 
we see that we have to prove 
\begin{eqnarray}  \label{eqn:ID} &&
A_p(A',\{-\hat \alpha\}) A_p(\{-\hat \beta\},B') 
(1-p^{-1+\hat\alpha+\hat\beta })
\sum_{j=0}^\infty \frac{\tau_{A'\cup\{-\hat\beta\}}(p^j)\tau_{B'\cup\{-\hat\alpha\}}(p^j)}{p^j}
\\
&& \qquad = \nonumber
\sum_{d,q=0}^\infty \frac{\mu(p^q)G_A(1-\hat \alpha, p^{d+q})G_B(1-\hat \beta,p^{d+q})}
{p^{d+q(2-\hat\alpha-\hat \beta)}}.
\end{eqnarray}
We note an easily proven identity that will be useful: if $a\in A$ and $A'=A-\{a\}$ then for $r\ge 1$
\begin{eqnarray} \label{eqn:tauid}
\tau_A(p^r)=\tau_{A'}(p^r)+p^{-a}\tau_A(p^{r-1}).
\end{eqnarray}
Thus, for  $r\ge 1$ we have 
\begin{eqnarray*}
G_A(s,p^{r})&=&
g_A(s,p^{r})\frac{p}{p-1}-\frac{p^{s}}{p-1}g_A(s ,p^{r-1})\\
&=& \prod_{\alpha\in A}( 1-p^{-s-\alpha})\frac{p}{p-1}\left( \sum_{j=0}^\infty \frac{\tau_A(p^{j+r})}{p^{js}}
-p^{s-1}\sum_{j=0}^\infty \frac{\tau_A(p^{j+r-1})}{p^{js}}\right) ;
\end{eqnarray*}
and with $s=1-a$ we have by (\ref{eqn:tauid})
\begin{eqnarray*}
G_A(1-a,p^{r})
 =\prod_{\alpha\in A'}( 1-p^{-1+a-\alpha} ) \sum_{j=0}^\infty \frac{\tau_{A'}(p^{j+r})}{p^{j(1-a)}}
 .
\end{eqnarray*}
The  right side of (\ref{eqn:ID}) is 
\begin{eqnarray} \label{eqn:goal}
\sum_{d=0}^\infty \frac{ G_A(1-\hat \alpha, p^{d})G_B(1-\hat \beta,p^{d})}
{p^{d}}-p^{-2+\hat \alpha+\hat \beta }\sum_{d=0}^\infty \frac{ G_A(1-\hat \alpha, p^{d+1})G_B(1-\hat \beta,p^{d+1})}
{p^{d}}
\end{eqnarray}
Now  
\begin{eqnarray}&& \label{eqn:sqdiff}
G_A(1-\hat \alpha, p^{d})G_B(1-\hat \beta,p^{d})
 -p^{-2+\hat\alpha+\hat\beta }   G_A(1-\hat \alpha, p^{d+1})G_B(1-\hat \beta,p^{d+1})\\ \nonumber
&&\qquad = \left(G_A(1-\hat \alpha, p^{d})-p^{-1+\hat \alpha} G_A(1-\hat \alpha, p^{d+1})\right)G_B(1-\hat \beta,p^{d})\\
&&\qquad \qquad + p^{-1+\hat\alpha} G_A(1-\hat \alpha, p^{d+1})\left(G_B(1-\hat \beta,p^{d})-p^{-1+\hat\beta}G_B(1-\hat \beta,p^{d+1})\right).
\nonumber
\end{eqnarray}
Using
\begin{eqnarray*}&&
G_A(1-\hat \alpha, p^{d})-p^{-1+\hat\alpha} G_A(1-\hat \alpha, p^{d+1})\\
&& \qquad =
\prod_{\alpha\in A'} (1-p^{-1+\hat\alpha-\alpha})\left(\sum_{j=0}^\infty \frac{\tau_{A'}(p^{j+d})}{p^{j(1-\hat \alpha)}}-p^{-1+\hat\alpha}
\sum_{j=0}^\infty \frac{\tau_{A'}(p^{j+d+1})}{p^{j(1-\hat\alpha)}}
\right)\\
&& \qquad = 
\prod_{\alpha\in A'} (1-p^{-1+a-\alpha})  \tau_{A'}(p^{d}) 
\end{eqnarray*}
we find that
\begin{eqnarray*}&&
\sum_{d=0}^\infty  \left(G_A(1-\hat \alpha, p^{d})-p^{-1+\hat \alpha} G_A(1-\hat \alpha, p^{d+1})\right)G_B(1-\hat \beta,p^{d})\\
&&\qquad = \prod_{\alpha\in A'} (1-p^{-1+\hat \alpha-\alpha})
\prod_{\beta\in B'} (1-p^{-1+\hat \beta-\beta})\sum_{d=0}^\infty \frac{\tau_{A'}(p^{d})}{p^d}\sum_{j=0}^\infty 
\frac{\tau_{B'}(p^{j+d})}{p^{j(1-\hat\beta)}}.
\end{eqnarray*}
The sum over $d$ and $j$ here may be rewritten as
\begin{eqnarray} \label{eqn:temp}
 \sum_{r=0}^\infty \frac{\tau_{B'}(p^r)}{p^{r }}
\sum_{d=0}^r p^{(r-d)\hat \beta }\tau_{A'}(p^d).
\end{eqnarray}
We recognize that this last sum over $d$ is a convolution:
\begin{eqnarray*} 
\sum_{d=0}^r p^{(r-d)\beta }\tau_{A'}(p^d)=\sum_{g\mid p^r} (p^r/g)^{-\hat \beta}\tau_{A'}(g)=\tau_{A'\cup \{-\hat\beta\}}(p^r).
\end{eqnarray*}
Thus, (\ref{eqn:temp}) is 
\begin{eqnarray*}
 \sum_{r=0}^\infty \frac{\tau_{B'}(p^r)\tau_{A'\cup \{-\hat\beta\}}(p^r)}{p^{r }}
 .
\end{eqnarray*}
The second term on the right side of (\ref{eqn:sqdiff}) is slightly different; it is
\begin{eqnarray*}
\prod_{\alpha\in A'} (1-p^{-1+\hat \alpha-\alpha})
\prod_{\beta\in B'} (1-p^{-1+\hat \beta-\beta})p^{-1+\hat\alpha}\sum_{d=0}^\infty \frac{\tau_{B'}(p^{d})}{p^d}\sum_{j=0}^\infty 
\frac{\tau_{A'}(p^{j+d+1})}{p^{j(1-\hat\alpha)}}.
\end{eqnarray*}
Now
\begin{eqnarray*}
  p^{-1+\hat\alpha}\sum_{d=0}^\infty \frac{\tau_{B'}(p^{d})}{p^d}\sum_{j=0}^\infty 
\frac{\tau_{A'}(p^{j+d+1})}{p^{j(1-\hat\alpha)}}&=&
\sum_{d=0}^\infty \frac{\tau_{B'}(p^{d})}{p^d}\sum_{j=1}^\infty 
\frac{\tau_{A'}(p^{j+d})}{p^{j(1-\hat\alpha)}}\\
&=& \sum_{d=0}^\infty \frac{\tau_{B'}(p^{d})}{p^d}\sum_{j=0}^\infty 
\frac{\tau_{A'}(p^{j+d})}{p^{j(1-\hat\alpha)}}-
\sum_{d=0}^\infty \frac{\tau_{B'}(p^{d}) \tau_{A'}(p^{d}) }{p^d} ;
\end{eqnarray*}
using the argument above this is 
\begin{eqnarray*}
\sum_{r=0}^\infty \frac{\tau_{A'}(p^r)\tau_{B'\cup \{-\hat\alpha\}}(p^r)}{p^{r }}
-\sum_{r=0}^\infty \frac{\tau_{A'}(p^r)\tau_{B' }(p^r)}{p^{r }}.
\end{eqnarray*}
Thus, (\ref{eqn:goal}) is
\begin{eqnarray*}
\sum_{r=0}^\infty \frac{\tau_{B'}(p^r)\tau_{A'\cup \{-\hat\beta\}}(p^r)}{p^{r }}
+\sum_{r=0}^\infty \frac{\tau_{A'}(p^r)\tau_{B'\cup \{-\hat\alpha\}}(p^r)}{p^{r }}
-\sum_{r=0}^\infty \frac{\tau_{A'}(p^r)\tau_{B' }(p^r)}{p^{r }}.
\end{eqnarray*}
Now 
\begin{eqnarray*}&&
 \tau_{B'}(p^r)\tau_{A'\cup \{-\hat\beta\}}(p^r) 
+ \tau_{A'}(p^r)\tau_{B'\cup \{-\hat\alpha\}}(p^r) 
- \tau_{A'}(p^r)\tau_{B' }(p^r) \\
&&\qquad =  \tau_{B'}(p^r)\left(\tau_{A'}(p^r)+p^{\hat\beta}\tau_{A'\cup \{-\hat\beta\}}(p^{r-1})\right)\\
&& \qquad \qquad \qquad 
+ \tau_{A'}(p^r)\left(\tau_{B'}(p^r)+p^{\hat\alpha}\tau_{B'\cup \{-\hat\alpha\}}(p^{r-1})\right)
- \tau_{A'}(p^r)\tau_{B' }(p^r) \\
&&\qquad = \tau_{A'}(p^r)\tau_{B' }(p^r) +
p^{\hat\beta}\tau_{B'}(p^r) \tau_{A'\cup \{-\hat\beta\}}(p^{r-1})
+p^{\hat\alpha}\tau_{A'}(p^r) \tau_{B'\cup \{-\hat\alpha\}}(p^{r-1})\\
&&\qquad = \left(\tau_{A'}(p^r)+p^{\hat\beta} \tau_{A'\cup \{-\hat\beta\}}(p^{r-1})\right)
\left(\tau_{B'}(p^r)+p^{\hat\alpha} \tau_{B'\cup \{-\hat\alpha\}}(p^{r-1})\right)
\\&& \qquad \qquad \qquad -p^{\hat\alpha+\hat\beta}
\tau_{A'\cup \{-\hat\beta\}}(p^{r-1}) \tau_{B'\cup \{-\hat\alpha\}}(p^{r-1})\\
&& \qquad = \tau_{A'\cup \{-\hat\beta\}}(p^{r}) \tau_{B'\cup \{-\hat\alpha\}}(p^{r})
-p^{\hat\alpha+\hat\beta}
\tau_{A'\cup \{-\hat\beta\}}(p^{r-1}) \tau_{B'\cup \{-\hat\alpha\}}(p^{r-1}).
\end{eqnarray*}
This leads to 
\begin{eqnarray*}
&&
1+\sum_{r=1}^\infty \frac{\tau_{A'\cup \{-\hat\beta\}}(p^{r}) \tau_{B'\cup \{-\hat\alpha\}}(p^{r})
-p^{\hat\alpha+\hat\beta}
\tau_{A'\cup \{-\hat\beta\}}(p^{r-1}) \tau_{B'\cup \{-\hat\alpha\}}(p^{r-1})}{p^r}\\&&\qquad 
= \left(1-p^{-1+\hat\alpha+\hat \beta}\right)\sum_{r=0}^\infty \frac{\tau_{A'\cup \{-\hat\beta\}}(p^{r}) \tau_{B'\cup \{-\hat\alpha\}}(p^{r})}{p^r};
\end{eqnarray*}
and (\ref{eqn:ID}) follows.

 \end{document}